# Minimum $L_1$-norm Estimation for Fractional Ornstein-Uhlenbeck Type Process Driven by a Hermite Process

B.L.S. Prakasa Rao [1]

CR Rao Advanced Institute of Mathematics, Statistics

and Computer Science, Hyderabad, India

**Abstract:** We investigate the asymptotic properties of the minimum $L_1$-norm estimator of the drift parameter for fractional Ornstein-Uhlenbeck type process driven by a Hermite process.

## 1 Introduction

We investigate the asymptotic properties of the minimum $L_1$-norm estimator of the drift parameter for a fractional Ornstein-Uhlenbeck type process driven by a Hermite process extending our earlier work in Prakasa Rao (2005) for processes driven by a fractional Brownian motion. Diffusion processes and diffusion type processes satisfying stochastic differential equations driven by Wiener processes are used for stochastic modeling in a wide variety of sciences such as population genetics, economic processes, signal processing as well as for modeling sunspot activity and more recently in mathematical finance. Statistical inference for diffusion type processes satisfying stochastic differential equations driven by Wiener processes have been studied earlier and a comprehensive survey of various methods is given in Prakasa Rao (1999). There has been an interest to study similar problems for stochastic processes, driven by a fractional Brownian motion, to model processes involving long range dependence (cf. Prakasa Rao (2010)). Le Breton (1998) studied parameter estimation and filtering in a simple linear model driven by a fractional Brownian motion. Kleptsyna and Le Breton (2002) studied parameter estimation problems for fractional Ornstein-Uhlenbeck process. The fractional Ornstein-Uhlenbeck process is a fractional analogue of the Ornstein-Uhlenbeck process, that is, a continuous time first order auto-regressive process $X = \{X_t, t \geq 0\}$ which is the solution of a one-dimensional homogeneous linear stochastic differential equation driven by a fractional Brownian motion (fBm) $W^H = \{W_t^H, t \geq 0\}$ with Hurst parameter $H$. Such a process is the unique Gaussian process satisfying the linear integral equation

(1. 1) $$X_t = x_0 + \theta \int_0^t X_s ds + \sigma W_t^H, t \geq 0.$$

They investigate the problem of estimation of the parameters $\theta$ and $\sigma^2$ based on the observation $\{X_s, 0 \leq s \leq T\}$ and study the asymptotic behaviour of these estimators as $T \to \infty$. Prakasa Rao (2024b) investigated the problem of parameter estimation for the fractional Ornstein-Uhlenbeck process when the driving force is a Gaussian process. Advances in statistical inference, for processes driven by fractional processes such as sub-fractional Brownian motion, fractional Levy process and mixed fractional Brownian motion, are surveyed in Prakasa Rao (2026).

In spite of the fact that maximum likelihood estimators (MLE) are consistent and asymptotically normal and also asymptotically efficient in general, they have some short comings at the same time. Their calculation is often cumbersome as the expression for MLE involve stochastic integrals at times which need good approximations for computational purposes. Furthermore MLE are not robust in the sense that a slight perturbation in the noise component will change the properties of MLE substantially. In order to circumvent such problems, the minimum distance approach is proposed. Properties of the minimum distance estimators (MDE) were discussed in Millar (1984) in a general frame work. Kutoyants and Pilibossian (1994) studied the problem of minimum $L_1$-norm estimation for the Ornstein-Uhlenbeck process. Prakasa Rao (2005) investigated the problem of minimum $L_1$-norm estimation for the fractional Ornstein-Uhlenbeck process driven by a fractional Brownian motion. El Machkouri et al. (2015), Chen and Zhou (2020), Lu (2022) study parameter estimation for an Ornstein-Uhlenbeck type process driven by a general Gaussian process. Nonparametric estimation of linear multiplier in stochastic differential equations driven by general Gaussian processes is studied in Prakasa Rao (2024a).

Hermite processes is another class of non-Gaussian processes which have been proposed as the driving forces for modeling stochastic phenomenon with long range dependence in a non-Gaussian environment (cf. Taqqu (1978); Lawrance and Kottegoda (1977)). Statistical inference for Vasicek-type model driven by Hermite process was studied in Nourdin and Diu Tran (2017) (cf. Diu Tran (2018)). Coupek and Kriz (2025) discuss parametric and nonparametric inference for nonlinear stochastic differential equations driven by Hermite processes. Our aim in this paper is to obtain the minimum $L_1$-norm estimator of the drift parameter of a fractional Ornstein-Uhlenbeck type process driven by a Hermite process and investigate the asymptotic properties of such estimators.

## 2 Hermite Processes

As the definition and properties of Hermite processes are not widely known, we now give a short review of the properties of Hermite processes for completeness following Diu Tran (2018) and Tudor (2013).

Let $W = \{W(h), h \in L^2(R)\}$ be a Brownian field defined on a probability space $(\Omega, \mathcal{F}, \mathcal{P})$, that is, a centered Gaussian family of random variables satisfying $E[W(h)W(g)] = <h, g>_{L^2(R)}$ for any $h, g \in L^2(R)$. For every $q \geq 1$, the $q$-th Wiener chaos $\mathcal{H}_q$ is defined as the closed linear subspace of $L^2(\Omega)$ generated by the family of random variables $\{H_q(W(h)), h \in L^2(R), \|h\|_{L^2(R)} = 1\}$, where $H_q$ is the $q$-th Hermite polynomial. Recall that $H_1(x) = x, H_2(x) = x^2 - 1, H_3(x) = x^3 - 3x, \ldots$. The mapping $I_q^W(h) = H_q(W(h))$ can be extended to a linear isometry between $L_s^2(R^q)$, the space of symmetric square integrable functions on $R^q$ equipped with the modified norm $\sqrt{q!} \|\cdot\|_{L^2(R^q)}$, and the $q$-th Wiener chaos $\mathcal{H}_q$. If $f \in L_s^2(R^q)$, , then the random variable $I_q^W(f)$, defined below, is called the *multiple Wiener integral* of the function $f$ of order $q$ and it is denoted by

$$I_q^W(f) = \int_{R^q} f(\psi_1, \ldots, \psi_q) dW_{\psi_1} \ldots dW_{\psi_q}.$$

For the existence and the properties of such multiple Wiener integrals, see Nualart (2006) and Major (2014).

We now define the Hermite process and discuss its properties following Diu Tran (2018) and Tudor (2013).

**Definition 2.1:** The *Hermite process* $\{Z_t^{q,H}, t \geq 0\}$ of order $q \geq 1$ and the self-similarity parameter $H \in (\frac{1}{2}, 1)$ is defined as

(2. 1) $$Z_t^{q,H} = c(q,H) \int_{R^q} \left( \int_0^t \Pi_{j=1}^q (s - \psi_j)_+^{H_0 - \frac{3}{2}} ds \right) dW_{\psi_1} \ldots dW_{\psi_q}$$

where

(2. 2) $$c(q,H) = \sqrt{\frac{H(2H-1)}{q!(\beta(H_0 - \frac{1}{2}, 2 - 2H_0))^q}}$$

and

(2. 3) $$H_0 = 1 + \frac{H-1}{q}.$$

Note that $H_0 \in (1 - \frac{1}{2q}, 1)$. The integral in (2.1) is a multiple Wiener integral as defined above of order $q$. Observe that $Z_0^{q,H} = 0$ a.s. The constant $c(q,H)$ is chosen in (2.1) so that $E[(Z_1^{q,H})^2] = 1$. Hermite process of order $q = 1$ is the fractional Brownian motion. It is the only Hermite process that is Gaussian. The Hermite process of order $q = 2$ is called the *Rosenblatt process*. The following properties are satisfied by the Hermite processes.

For any two processes $X = \{X_t, t \geq 0\}$ and $Y = \{Y_t, t \geq 0\}$, we say that $\{X_t, t \geq 0\} \triangleq \{Y_t, t \geq 0\}$ if the finite dimensional distributions of the process $X$ are the same as the finite dimensional distributions of the process $Y$.

**Theorem 2.1:** Suppose $Z = \{Z_t^{q,H}, t \geq 0\}$ is the Hermite process of order $q \geq 1$ and the parameter $H \in (\frac{1}{2}, 1)$. Then the process is centered and the following properties hold:

(i) the process $Z$ is self-similar, that is, for all $a > 0$

$$\{Z_{at}^{q,H}, t \geq 0\} \triangleq \{a^H Z_t^{q,H}, t \geq 0\};$$

(ii) the process $Z$ has stationary increments, that is, for any $h > 0$,

$$\{Z_{t+h}^{q,H} - Z_h^{q,H}, t \geq 0\} \triangleq \{Z_t^{q,H}, t \geq 0\};$$

(iii) the covariance function of the process $Z$ is given by

$$E[Z_t^{q,H} Z_s^{q,H}] = \frac{1}{2}(t^{2H} + s^{2H} - |t-s|^{2H}), t \geq 0, s \geq 0;$$

(iv) the process $Z$ has long range dependence property, that is,

$$\sum_{n=0}^{\infty} E[Z_1^{q,H}(Z_{n+1}^{q,H} - Z_n^{q,H}) = \infty;$$

(v) for any $\alpha \in (0, H)$ and any interval $[0,T] \subset R_+$, the process $\{Z_t^{q,H}, 0 \leq t \leq T\}$ admits a version with Holder continuous sample paths of order $\alpha$; and

(vi) for every $p \geq 1$, there exists a positive constant $C_{p,q,H}$ such that

$$E[|Z_t^{q,H}|^p] \leq C_{p,q,H} t^{pH}, t \geq 0.$$

For proof of this theorem, see Tudor (2013).

**Wiener integral with respect to the Hermite process:**

The Wiener integral of a non-random function $f$, with respect to the Hermite process $Z = \{Z_t^{q,H}, t \geq 0\}$, is studied in Maejima and Tudor (2007) and it is denoted by

$$E[\int_R f(u)dZ_u^{q,H}$$

whenever it exists. It is known that the Wiener integral of the function $f$ with respect to the Hermite process of order $q$ and parameter $H$ exists if

(2. 4) $$\int_R \int_R |f(u)f(v)||u-v|^{2H-2} \, dudv < \infty.$$

Let $|\mathcal{H}|$ denote the class of functions $f$ satisfying the inequality (2.4). It can be shown that , for any $f, g \in |\mathcal{H}|$,

(2. 5) $$E[\int_R f(u)dZ_u^{q,H} \int_R g(u)dZ_u^{q,H}] = H(2H-1) \int_R \int_R f(u)g(v)|u-v|^{2H-2}dudv.$$

Furthermore, the Wiener integral of $f$ with respect to the process $Z$ admits the representation

(2. 6) $$\int_R f(u)dZ_u^{q,H} = c(q,H) \int_{R^q} (\int_R f(u)\Pi_{j=1}^q (u-\psi_j)_+^{H_0-\frac{3}{2}} du)dW_{\psi_1} \ldots dW_{\psi_q}$$

where $c(q,H)$ and $H_0$ are as defined by the equations (2.2) and (2.3) respectively.

## 3 Minimum $L_1$-norm Estimation

Let $(\Omega, \mathcal{F}, (\mathcal{F}_t), P)$ be a stochastic basis satisfying the usual conditions and the processes discussed in the following are $(\mathcal{F}_t)$-adapted. Further the natural filtration of a process is understood as the $P$-completion of the filtration generated by this process. Let $Z = \{Z_t^{q,H}, t \geq 0\}$ denote the Hermite process of order $q \geq 1$ and self-similarity parameter $H \in (\frac{1}{2}, 1)$. This process is centered, $H$-self-similar, has stationary increments and exhibits long-range dependence. If $q = 1$, then the process is the fractional Brownian motion with Hurst index $H$ and if $q = 2$, then the process is called a Rosenblatt process. If $q \geq 2$, then the process is non-Gaussian.

Let us consider a stochastic process $\{X_t, t \in [0,1]\}$ defined by the stochastic integral equation

(3. 1) $$X_t = \psi + \theta \int_0^t X_s ds + \varepsilon Z_t^{q,H}, 0 \leq t \leq 1,$$

on the probability space $(\Omega, \mathcal{F}, (\mathcal{F}_t), P)$ where $\theta$ is an unknown drift parameter, $\varepsilon > 0$ and $\psi$ is a random variable independent of the process $X$. Let $P_\theta^{(\varepsilon)}$ be the probability measure generated by the process $\{X_t, t \in [0,1]\}$ when $\theta$ is the true parameter. For convenience, we write the above integral equation in the form of a stochastic differential equation

(3. 2) $$dX_t = \theta X_t dt + \varepsilon dZ_t^{q,H}, X_0 = \psi, 0 \leq t \leq 1, (3)$$

driven by the Hermite process $Z$. There exists a unique continuous solution of the equation (3.2) by the results in Maejima and Tudor (2007). Existence and uniqueness of the solution $\{X_t, 0 \leq t \leq 1\}$, for the differential equation given above, also follows from Theorem 5.1 in Coupek and Kriz (2025) under some conditions. The unique continuous solution of (3.2) is given by

$$X_t = e^{\theta t}(\psi + \int_0^t e^{-\theta u} dZ_u^{q,H}), t \geq 0.$$

In particular, if we choose $\psi = \varepsilon \int_{-\infty}^t e^{-\theta u} dZ_u^{q,H}$, then the unique continuous solution of (3.2) is

(3. 3) $$X_t = \varepsilon \int_{-\infty}^t e^{\theta(t-u)} dZ_u^{q,H}, t \geq 0.$$

If the initial value $\psi$ is taken to be zero, then the unique continuous solution is given by

(3. 4) $$X_t = \varepsilon \int_0^t e^{\theta(t-u)} dZ_u^{q,H}, t \geq 0.$$

We now consider the problem of estimation of the parameter $\theta$ based on the observation of the Hermite fractional Ornstein-Uhlenbeck type process $X = \{X_t, 0 \leq t \leq 1\}$ satisfying the stochastic differential equation

(3. 5) $$dX_t = \theta X_t dt + \varepsilon dZ_t^{q,H}, X_0 = x_0, 0 \leq t \leq 1$$

where $\theta \in \Theta \subset R, \varepsilon > 0, \{Z_t^{q,H}, 0 \leq t \leq 1\}$ is the Hermite process of order $q \geq 1$ and the self-similarity parameter $H \in (\frac{1}{2}, 1)$ and study the asymptotic properties of the estimator as $\varepsilon \to 0$. We assume that the parameters $q$ and $H$ and the initial values $x_0$ are known.

Let $x_t(\theta)$ be the solution of the corresponding ordinary differential equation

$$\frac{dx_t(\theta)}{dt} = \theta x_t$$

with $x_t(\theta) = x_0$ for $t = 0$. It is obvious that

(3. 6) $$x_t(\theta) = x_0 e^{\theta t}, 0 \leq t \leq 1.$$

Let

(3. 7) $$S(\theta) = \int_0^1 |X_t - x_t(\theta)| dt.$$

We define $\hat{\theta}_\varepsilon$ to be a *minimum $L_1$-norm estimator* if there exists a measurable selection $\hat{\theta}_\varepsilon$ such that

(3. 8) $$S(\hat{\theta}_\varepsilon) = \inf_{\theta \in \Theta} S(\theta).$$

Conditions for the existence of a measurable selection are given in Lemma 3.1.2 in Prakasa Rao (1987). We assume that there exists a measurable selection $\hat{\theta}_\varepsilon$ satisfying the above equation. An alternate way of defining the estimator $\hat{\theta}_\varepsilon$ is by the relation

(3. 9) $$\hat{\theta}_\varepsilon = \arg\inf_{\theta \in \Theta} \int_0^1 |X_t - x_t(\theta)| dt.$$

Let $Z_t^{q,H,*} = \sup_{0 \leq s \leq t} |Z_s^{q,H}|$. Maximal inequalities for processes such as the fractional Brownian motion and sub-fractional Bownian motion are known and are reviewed in Prakasa Rao (2014, 2017, 2020) following the property of self-similarity for such processes. Maximal inequalities for general Gaussian processes are presented in Li and Shao (2001), Berman (1985), Marcus and Rosen (2006) and Borovkov et al. (2017) among others. As far as we are aware, there are no maximal inequalities available for Hermite processes in the literature. We now derive a maximal inequality for the Hermite process as a consequence of the self-similarity property of the process.

For any process $X = \{X_t, t \geq 0\}$, define
$$X_t^* = \sup\{X_s, 0 \leq s \leq t\}.$$

Since the Hermite process $\{Z_t^{q,H}, t \geq 0\}$ is a self-similar process with self-similarity parameter $H$, it follows that
$$Z_{at}^{q,H} \stackrel{\Delta}{=} a^H Z_t^{q,H}$$
for every $t \geq 0$ and $a > 0$. This in turn implies that
$$Z_{at}^{q,H,*} \stackrel{\Delta}{=} a^H Z_t^{q,H,*}$$
for every $a > 0$. Hence the following result holds.

**Theorem 3.1:** Let $T > 0$ and $\{Z_t^{q,H}, t \geq 0\}$ be the Hermite process with self-similarity parameter $H$. Suppose that $E[|Z_1^{q,H,*}|^p] < \infty$ for some $p > 0$.. Then
$$E[(Z_T^{q,H,*})^p] = K(p,q,H) T^{pH}$$
where $K(p,q,H) = E[(Z_1^{q,H,*})^p]$.

**Proof:** The result is a consequence of the fact $Z_T^{q,H,*} = T^H Z_1^{q,H,*}$ by self-similarity.

# 4 Consistency of the estimator:

Let $\theta_0$ denote the true parameter. For any $\delta > 0$, define

(4. 1) $$g(\delta) = \inf_{|\theta - \theta_0| > \delta} \int_0^1 |x_t(\theta) - x_t(\theta_0)| dt.$$

Note that $g(\delta) > 0$ for any $\delta > 0$.

**Theorem 4.1:** Suppose $\{Z_t^{q,H}, t \geq 0\}$ is a Hermite process on the interval $[0,1]$ such that $E|Z_1^{q,H,*}| < \infty$. Let $m = E[\sup_{0 \leq t \leq 1} |Z_t^{q,H}|]$. Then there exists a positive constant $C$ depending on $q, H$ such that for every $\delta > 0$,

$$P_{\theta_0}^{(\varepsilon)}(|\hat{\theta}_\varepsilon - \theta_0| > \delta) \leq C\varepsilon.$$

**Proof:** Let $||.||$ denote the $L_1$-norm. Then

(4. 2) $$\begin{aligned} P_{\theta_0}^{(\varepsilon)}(|\hat{\theta}_\varepsilon - \theta_0| > \delta) &= P_{\theta_0}^{(\varepsilon)}(\inf_{|\theta - \theta_0| \leq \delta} ||X - x(\theta)|| > \inf_{|\theta - \theta_0| > \delta} ||X - x(\theta)||) \\ &\leq P_{\theta_0}^{(\varepsilon)}(\inf_{|\theta - \theta_0| \leq \delta}(||X - x(\theta_0)|| + ||x(\theta) - x(\theta_0)||) \\ &\quad > \inf_{|\theta - \theta_0| > \delta}(||x(\theta) - x(\theta_0)|| - ||X - x(\theta_0)||)) \\ &= P_{\theta_0}^{(\varepsilon)}(2||X - x(\theta_0)|| > \inf_{|\theta - \theta_0| > \delta} ||x(\theta) - x(\theta_0)||) \\ &= P_{\theta_0}^{(\varepsilon)}(||X - x(\theta_0)|| > \frac{1}{2} g(\delta)). \end{aligned}$$

Since the process $X$ satisfies the stochastic differential equation (3.2), it follows that

(4. 3) $$\begin{aligned} X_t - x_t(\theta_0) &= x_0 + \theta_0 \int_0^t X_s ds + \varepsilon Z_t^{q,H} - x_t(\theta_0) \\ &= \theta_0 \int_0^t (X_s - x_s(\theta_0)) ds + \varepsilon Z_t^{q,H} \end{aligned}$$

a.s. $P_{\theta_0}^{(\varepsilon)}$ since $x_t(\theta) = x_0 e^{\theta t}$. Let $U_t = X_t - x_t(\theta_0)$. Then it follows from the equation given above that

(4. 4) $$U_t = \theta_0 \int_0^t U_s \, ds + \varepsilon Z_t^{q,H}.$$

Let $V_t = |U_t| = |X_t - x_t(\theta_0)|$. The relation given above implies that

(4. 5) $$V_t = |X_t - x_t(\theta_0)| \leq |\theta_0| \int_0^t V_s ds + \varepsilon |Z_t^{q,H}|.$$

Applying the Gronwall-Bellman Lemma, it follows that

(4. 6)
$$\sup_{0\leq t\leq 1} |V_t| \leq \varepsilon e^{|\theta_0|} \sup_{0\leq t\leq 1} |Z_t^{q,H}|.$$

Hence

(4. 7)
$$P_{\theta_0}^{(\varepsilon)}(\|X - x(\theta_0)\| > \frac{1}{2}g(\delta)) \leq P[\sup_{0\leq t\leq 1} |Z_t^{q,H}| > \frac{e^{-|\theta_0|}g(\delta)}{2\varepsilon}]$$
$$= P[Z_1^{q,H,*} > \frac{e^{-|\theta_0|}g(\delta)}{2\varepsilon}].$$

Let $m = E[\sup_{0\leq t\leq 1} |Z_t^{q,H}|]$. From the maximal inequality proved in Theorem 3.1 for Hermite processes, it follows that

$$m \leq C(q, H)$$

for some positive constant $C(q, H)$. For fixed $\delta > 0$, we can choose $\epsilon$ sufficiently small so that $\frac{e^{-|\theta_0|}g(\delta)}{2\varepsilon} > m$. For such $\varepsilon$,

(4. 8)
$$P_{\theta_0}^{(\varepsilon)}[|\hat{\theta}_\varepsilon - \theta_0| > \delta] \leq \frac{2\varepsilon}{e^{-|\theta_0|}g(\delta)} E[Z_1^{q,H,*}]$$
$$\leq C_1(\delta, q, H)\epsilon$$

for some positive constant $C_1(\delta, q, H)$ independent of $\varepsilon$.

**Remarks:** As a consequence of the result obtained above, it follows that

$$P_{\theta_0}^{(\varepsilon)}(|\hat{\theta}_\varepsilon - \theta_0| > \delta) \to 0 \text{ as } \epsilon \to 0$$

for every $\delta > 0$. Hence the minimum norm $L_1$-estimator $\hat{\theta}_\varepsilon$ is weakly consistent for estimating the parameter $\theta_0$.

## 5 Asymptotic distribution of the estimator

We will now study the asymptotic distribution if any of the estimator $\hat{\theta}_\varepsilon$ after suitable scaling. It can be checked that

(5. 1)
$$X_t = e^{\theta_0 t}\{x_0 + \int_0^t e^{-\theta_0 s}\varepsilon dZ_s^{q,H}\}$$

or equivalently

(5. 2)
$$X_t - x_t(\theta_0) = \varepsilon e^{\theta_0 t} \int_0^t e^{-\theta_0 s}dZ_s^{q,H}.$$

Let
(5. 3) $$Y_t = e^{\theta_0 t}\int_0^t e^{-\theta_0 s}dZ_s^{q,H}.$$

The integral with respect to the process $\{Z_t^{q,H}, t \geq 0\}$ is interpreted as a Wiener integral (cf. Maejima and Tudor (2007); Coupek and Kriz (2025)). For any $0 \leq t, s \leq 1$,

(5. 4) $$\begin{aligned}Cov(Y_t, Y_s) &= e^{\theta_0 t + \theta_0 s}E[\int_0^t e^{-\theta_0 u}dZ_u^{q,H}\int_0^s e^{-\theta_0 v}dZ_v^{q,H}]\\ &= R(t,s) \text{ (say)}.\end{aligned}$$

In particular
(5. 5) $$Var(Y_t) = R(t,t).$$

Observe that $\{Y_t, 0 \leq t \leq 1\}$ is a zero mean process with $Cov(Y_t, Y_s) = R(t,s)$. Let

(5. 6) $$\zeta = \arg\inf_{-\infty < u < \infty}\int_0^1 |Y_t - utx_0 e^{\theta_0 t}|dt.$$

**Theorem 5.1:** As $\varepsilon \to 0$, the random variable $\varepsilon^{-1}(\hat{\theta}_\varepsilon - \theta_0)$ converges in probability to a random variable whose probability distribution is the same as that of the random variable $\zeta$ defined by the equation (5.6).

**Proof:** Let $x_t'(\theta) = x_0 t e^{\theta t}$ and let

(5. 7) $$J_\varepsilon(u) = ||Y - \varepsilon^{-1}(x(\theta_0 + \varepsilon u) - x(\theta_0))||$$

and
(5. 8) $$J_0(u) = ||Y - ux'(\theta_0)||.$$

Furthermore, let
(5. 9) $$A_\varepsilon = \{\omega : |\hat{\theta}_\varepsilon - \theta_0| < \delta_\varepsilon\}, \delta_\varepsilon = \varepsilon^\tau, \tau \in (\frac{1}{2}, 1), L_\varepsilon = \varepsilon^{\tau-1}.$$

Observe that the random variable $u_\varepsilon^* = \varepsilon^{-1}(\theta_\varepsilon^* - \theta_0)$ satisfies the equation

(5. 10) $$J_\varepsilon(u_\varepsilon^*) = \inf_{|u|<L_\varepsilon} J_\varepsilon(u), \omega \in A_\varepsilon.$$

Define
(5. 11) $$\zeta_\varepsilon = \arg\inf_{|u|<L_\varepsilon} J_0(u).$$

Furthermore note that, with probability one,

$$
\begin{aligned}
(5.12)\quad \sup_{|u|<L_\varepsilon} |J_\varepsilon(u) - J_0(u)| &= \sup_{|u|<L_\varepsilon} |||Y - ux'(\theta_0) - \frac{1}{2}\varepsilon u^2 x''(\tilde\theta)|| - ||Y - ux'(\theta_0)|||\\
&\leq \frac{\varepsilon}{2}L_\varepsilon^2 \sup_{|\theta-\theta_0|<\delta_\varepsilon} \int_0^1 |x''(t)|dt\\
&\leq C\varepsilon^{2\tau-1}.
\end{aligned}
$$

Here $\tilde\theta = \theta_0 + \alpha(\theta - \theta_0)$ for some $\alpha \in (0,1)$. Note that the last term in the above inequality tends to zero as $\varepsilon \to 0$. Furthermore the process $\{J_0(u), -\infty < u < \infty\}$ has a unique minimum $u^*$ with probability one. In addition, we can choose the interval $[-L, L]$ such that

$$(5.13)\qquad P^{(\varepsilon)}_{\theta_0}\{u^*_\varepsilon \in (-L, L)\} \geq 1 - \beta(g(L))^{-1}$$

and

$$(5.14)\qquad P\{u^* \in (-L, L)\} \geq 1 - \beta(g(L))^{-1}$$

where $\beta > 0$. Note that $g(L)$ increases as $L$ increases. The processes $\{J_\varepsilon(u), u \in [-L, L]\}$ and $J_0(u), u \in [-L, L]$ satisfy the Lipschitz conditions and $J_\varepsilon(u)$ converges uniformly to $J_0(u)$ over $u \in [-L, L]$. Hence the minimizer of $J_\varepsilon(.)$ converges to the minimizer of $J_0(u)$. This completes the proof.

**Remarks :** We have seen earlier that the process $\{Y_t, 0 \leq t \leq 1\}$ is a zero mean process with the covariance function $Cov(Y_t, Y_s) = R(t,s)$ for $0 \leq t, s \leq 1$. Recall that

$$(5.15)\qquad \zeta = \arg\inf_{-\infty<u<\infty} \int_0^1 |Y_t - utx_0 e^{\theta_0 t}|dt.$$

It is not clear what the distribution of the random variable $\zeta$ is. It depends on the Hermite process $Z^{q,H}$. Observe that for every $u$, the integrand in the above integral is the absolute value of the process $\{J_t, 0 \leq t \leq 1\}$ with the mean function $E(J_t) = -utx_0 e^{\theta_0 t}$ and the covariance function $Cov(J_t, J_s) = R(t,s)$ for $0 \leq s, t \leq 1$. It is easy to extend the results to the Hermite process defined on any interval $[0, T]$ for any fixed $T > 0$. The distribution of

$$u_T = \arginf_{u\in R}\int_0^T |Y_t - utx_0 e^{\theta_0 t}|dt,$$

as $T \to \infty$, has been investigated in Aubry (1999) for a diffusion process, Diop and Yode (2010) for the Ornstein-Uhlenbeck process driven by a Levy process and Kutoyants and Pilibossian (1994) for the Ornstein-Uhlenbeck process. Diu Tran (2018) proved that, if either $q \geq 2$ or $q = 1$ and $H \geq \frac{3}{4}$, then

$$T^{\frac{2}{q}(1-H)-1}\int_0^T (Y_t^2 - E(Y_t^2))dt \xrightarrow{L} B_{H,q}(-\theta_0)^{-2H-\frac{2}{q}(1-H)}\zeta$$

as $T \to \infty$ where $\zeta$ is distributed according to the Rosenblatt distribution with parameter $4H' = 1 - \frac{2}{q}(1-H)$ and $B_{H,q}$ is defined by

$$B_{H,q} = \frac{H(2H-1)}{\sqrt{(H_0 - \frac{1}{2})(4H_0 - 3)}} \frac{\Gamma(2H + \frac{2}{q}(1-H))}{2H + \frac{2}{q}(1-H) - 1}$$

and $H_0 = 1 - \frac{1-H}{q}$.

**Acknowledgment:** This work was supported by the INSA Honorary Scientist fellowship at the CR Rao Advanced Institute of Mathematics, Statistics and Computer Science, Hyderabad, India.

CR Rao Advanced Institute of Mathematics, Statistics and Computer Science, Hyderabad, India.
e-mail: blsprao@gmail.com